\documentclass[a4paper,fleqn]{cas-dc}

\usepackage[numbers]{natbib}

\def\tsc#1{\csdef{#1}{\textsc{\lowercase{#1}}\xspace}}
\tsc{WGM}
\tsc{QE}

\newdefinition{rmk}{Remark}
\newproof{pf}{Proof}
\newdefinition{definition}{Definition}
\newdefinition{example}{Example}

\shorttitle{Fractional Liu uncertain differential equation and its
	application to finance}    
\shortauthors{Alireza Najafi, Rahman Taleghani}  
\title [mode = title]{Comments for ‘‘Comments for Fractional Liu uncertain differential equation and its
	application to finance’’}  
\author[1]{Alireza Najafi}

\cormark[*]
\fnmark[*]
\ead{alireza.najafi@guilan.ac.ir}

\affiliation[1]{organization={Department of Applied Mathematics, Faculty of Mathematical Sciences, University of Guilan},
            city={Rasht},
            postcode={41938-1914}, 
            state={Guilan},
            country={Iran}}

\author[1]{Rahman Taleghani}

\ead{Rahman.taleghani@gmail.com}

\cortext[1]{Corresponding author}
\begin{document}

\begin{abstract}
In this article, we  demonstrate that the claim made by Xiaoyan Li and Ni Sun \cite{bib 1} regarding the incorrectness of Theorem 7 in the paper \cite{bib 2}  is wrong, and  show that this  Theorem is  based on the integral with respect to $(dt)^{\alpha}$  and the fractional derivative of order $\alpha$ ($ 0 < \alpha \leq 1$). 
\end{abstract}
\maketitle



	\section{Introduction}
	In 2005,  Guy Jumarie studied on the fractional Brownian motion and proposed a  representation of fractional Brownian motion as an integral with respect to $ (dt)^{\alpha}$. To do this, he used the fractional derivative of order $\alpha$ concepts \cite{bib 3}. After that, he  discussed about the  lagrange analytical mechanics and used the fractional derivative and the  integral with respect to $ (dt)^{\alpha}$ and   proposed an extension of the Lagrange analytical mechanics to deal with dynamics of fractal nature \cite{bib 4}. In 2010, Lv Longjin et al. studied about the  connection between fractional derivative and fractional Brownian motion  and used the  Jumarie's idea and presented an application of fractional derivatives in stochastic models driven by fractional Brownian motion \cite{bib 5}. In recent years, researchers used  the Jumarie's idea to find the option price PDE under the fractional  and mixed fractional financial models \cite{bib 6, bib 7, bib 8,bib 9}.

	In paper \cite{bib 2}, we used the integral with respect to $ (dt)^{\alpha}$  to show that the fractional uncertain differential equation has a unique solution.   Xiaoyan Li and Ni Sun in \cite{bib 1} by using the wrong definition of the  integral with respect to $ (dt)^{\alpha}$ proved that the Jumarie's idea and therefore  Theorem  7 that is  presented in \cite{bib 2} is wrong. Here, we study the Jumarie's idea and  a part of the   Theorem  7 that we used the idea.

	\section{Main result}
	Here, we  present the integral with respect to $ (dt)^{\alpha}$ which is introduced by  Guy Jumarie and after that we mention some examples to better understand this definition.  In 2005,  Guy Jumarie created a bridge between the  integral with respect to $ (dt)^{\alpha}$ and the fractional integral as follows : 
	\begin{definition}\label{def 1}
		Let $f : \mathbb{R} \longrightarrow \mathbb{R}$ be a continuous function. Then its integral w.r.t $(dt)^{\alpha}$ for $0 < \alpha \leq 1$ is defined as follows:
		\begin{align}\label{eq 11}
			\int_0^t f(\tau) (d\tau)^{\alpha} = \alpha \int_0^t (t - \tau)^{\alpha -1}f(\tau)  d\tau
		\end{align}
	\end{definition}
	This definition is based on the  fractional derivative of order $\alpha$ for $0 < \alpha \leq 1$ \cite{bib 3}. Also, to better understand this definition, he gave some examples as follows \cite{bib 3, bib 4}:
	\begin{example}
		.
		\\
		(i) $~~$ Let $f(\tau) = 1$, then 
		\begin{align}
			\int_0^t  (d\tau)^{\alpha} = t^{\alpha} ~~~~ 0 < \alpha \leq 1.
		\end{align}
		(ii) $~~$ Let $f(\tau) =\tau^{\gamma}$, then 
		\begin{align}
			\int_0^t \tau^{\gamma} (d\tau)^{\alpha} = \frac{\Gamma(\alpha + 1) ~ \Gamma(\gamma + 1)}{\Gamma(\alpha + \gamma + 1)} t^{\alpha + \gamma}, ~~~~ 0 < \alpha \leq 1.
		\end{align}
		\\
		(iii) $~~$ Let $f(t)$ be  the Dirac delta generalized function $\delta (t)$, then 
		\begin{align}
			\int_0^t \delta (\tau) (d\tau)^{\alpha} = \alpha t^{\alpha -1}, ~~~~ 0 < \alpha \leq 1.
		\end{align}
	\end{example}
	It should be noted that this definition is based on fractional calculations and it cannot be expressed as a Riemann sum like the Riemann integral or Riemann Stieltjes Integral. 
	\\
	In paper \cite{bib 2}, we studied  the  fractional uncertain differential equation and proved that the equation has a unique solution. In fact, let $(\Gamma, \mathcal{M})$ be an uncertain space and  $\left\lbrace X_t \right\rbrace_{t \geq 0}$ satisfies in  
	\begin{align}\label{eq 43}
		dX_t = g(t, X_t) dt + h(t, X_t) dF_t^{\alpha},
	\end{align}  
	where $\left\lbrace F_t^{\alpha} \right\rbrace_{t \geq 0}$ is the fractional Liu process and   the coefficients $g(t, x)$ and $h(t, x)$ satisfy the linear growth and Lipschitz  conditions. Then the Eq \eqref{eq 43} has a unique solution. Here we present a part of the proof of theorem 7.  Assume that  $X_t^{(0)} = X_0$, and 
	\begin{align}
		X_t^{(n)} = X_0 &+ \int_0^t g\left( s, X_s^{(n-1)}\right) ds + \int_0^t h\left( s, X_s^{(n-1)}\right) dF_s^{\alpha},
	\end{align}
	for $n =1, 2, \cdots$, and suppose 
	\begin{align}
		D_t^{(n)}(\gamma) = \underset{0 < s < t}{max} \vert X_s^{(n+1)} (\gamma) - X_s^{(n)} (\gamma) \vert, 
	\end{align}
	for each $\gamma \in \Gamma$. It follows from the linear growth  and Lipschitz conditions and the fractional Liu process property, we have
	\begin{align}
	&	D_t^{(0)}(\gamma)  = \underset{0 < s < t}{max} \vert \int_0^s g\left(u, X_0\right) du + \int_0^s h\left(u, X_0\right) dF_u^{\alpha} \vert \nonumber \\
		& <  \int_0^t \vert g\left(u, X_0\right) \vert  du  +{\color{red} \frac{2 K(\gamma)}{1 - \alpha}\int_0^t \vert h\left(u, X_0\right) \vert  (du)^{1-\alpha}} \nonumber \\ 
		& = \int_0^t \vert g\left(u, X_0\right) \vert  du \nonumber \\ 
		&  + {\color{red}\frac{2 K(\gamma)}{1 - \alpha}(1-\alpha)\int_0^t (t - u)^{-\alpha} \vert h\left(u, X_0\right) \vert du} \nonumber \\
		& < L(1+|X_0|) \left(t + \frac{2 t^{1-\alpha} K(\gamma)}{1-\alpha}\right),
	\end{align}
	In this relation, we used Definition \ref{def 1} and  deduced that 
	\begin{align}\label{eq 26}
		\frac{2 K(\gamma)}{1 - \alpha}\int_0^t& \vert h\left(u, X_0\right) \vert   (du)^{1-\alpha} = \nonumber \\
		& \frac{2 K(\gamma)}{1 - \alpha}(1-\alpha)\int_0^t (t - u)^{-\alpha} \vert h\left(u, X_0\right) \vert du, 
	\end{align}
	where $h$ is a continuous function. \\
	In \cite{bib 1}, researchers used the Riemann sum and showed that the equation \eqref{eq 11} and thus the equation  \eqref{eq 26} are wrong. In fact, they used the Riemann sum concepts and considered the integral with respect to $(dt)^{\alpha}$ as follows 
	\begin{align}
		\int_0^t f(s) (ds)^{\alpha} = \underset{n \longrightarrow \infty}{\lim} \sum_{i=1}^n f(s_i) (\Delta s_i)^{\alpha},
	\end{align}
	where $0 = s_0 < s_1 < \cdots < s_n = t$. Then for $\alpha = \frac{1}{2}$, $t = 1$,  and $f(s) = 1$ for any $s \in [0, 1]$,  they concluded
	\begin{align}
		\int_0^1  (ds)^{\alpha} =  \underset{n \longrightarrow \infty}{\lim} \sum_{i=1}^n (\Delta s_i)^{\alpha}	= \underset{n \longrightarrow \infty}{\lim} \sum_{i=1}^n \left(\frac{1}{n}\right)^{\alpha} \longrightarrow  \infty
	\end{align}
	and $0 = s_0 < s_1 < \cdots < s_n = 1$. Therefore,  the $\int_0^1  (ds)^{\alpha}$ is diverge, and thus  the definition \ref{def 1}  is meaningless.

	Based on this argument, they  concluded that Theorem 7 in \cite{bib 2} is wrong and discussed another equation and considered it as a replacement for this equation.

	But according to the definition \ref{def 1}, we have 
	\begin{align}
		\int_0^1  (ds)^{\alpha} = \alpha \int_0^1  (1 - s)^{\alpha -1} ds = 1. 
	\end{align}

	

\end{document}